\documentclass{amsart}
\usepackage{amsmath,amscd}
\usepackage{amssymb}
\usepackage{bbm}
\usepackage{enumerate}

\newtheorem{theorem}{Theorem}[section] 

\newtheorem{corollary}[theorem]{Corollary}

\theoremstyle{definition}
\newtheorem{definition}[theorem]{Definition}

\newcommand{\f}{\mathfrak}
\newcommand{\pro}{{\rm Pro}}
\newcommand{\mc}{\mathcal}
\theoremstyle{remark}

\numberwithin{equation}{section}

\begin{document}

\title{Projection of Elliptic Orbits and Branching Laws}

\author{Hongyu He}
\address{Department of Mathematics, Louisiana State University, Baton Rouge, Louisiana 70803}
\email{livingstone@alum.mit.edu}

%    General info

\date{}

\footnote{This article is based on the talk I gave at the conference on Representations and Characters: Revisiting the Works of Harish-Chandra and André Weil at the National University of Singapore in 2022.}

\begin{abstract}
Let $G$ be a Lie group, and $H\subset G$ a closed Lie subgroup. Let $\pi$ be an irreducible unitary representation of $G$. In this paper, we briefly discuss the orbit method and its application to the branching problem $\pi|_{H}$. We use the Gan-Gross-Prasad branching law for $(G, H)= ( U(p,q), U(p, q-1) )$ as an example to illustrate the relation between $\pro_{\f u(p, q-1)}^{\f u(p,q)} \mc O(\lambda)$ and the branching law of the discrete series $D_{\lambda}|_{U(p,q-1)}$ for $\lambda$ a regular elliptic element. We also discuss some results regarding branching laws and wave front sets.
The presentation of this paper does not follow the historical timeline of development.
\end{abstract}

\maketitle

\section{Coadjoint Orbits and Their Projections}
Let $G$ be a connected Lie group. Let $\f g$ be its Lie algebra, and $\f g^*$ its dual space (over $\mathbb R$). The Lie group $G$ acts on $\f g$ and $\f g^*$ respectively. Each orbit is called an adjoint orbit or a coadjoint orbit respectively. If $\mathfrak g$ is reductive, then $\f g$ can be identified with $\mathfrak g^*$ such that the adjoint action coincides with coadjoint action. Then adjoint orbits can be identified with coadjoint orbits. We denote the set of adjoint orbits by $\f g//G$ and the set of coadjoint orbits by $\f g^*//G$. For each $\lambda \in \f g^*$, let $\mc O(\lambda)$ be the corresponding orbit generated by $\lambda$. Since our discussion will not involve the topology on $\f g//G$ or $\f g^*//G$, we will not address  this issue in this paper. \\
\\
One classical example is the unitary Lie algebra $\f u(n)$, consisting of the $n \times n$ skew-Hermitian matrices.  The unitary group $U(n)$ acts on $\mathfrak u(n)$ adjointly. The adjoint orbits in this case are in one-to-one correspondence with $i \lambda$ with $\lambda \in \mathbb R^n$ and
$$\lambda_1 \geq \lambda_2 \ldots \geq \lambda_n.$$
\\
Let $H$ be a Lie subgroup of $G$. Let $\mathcal O(\lambda)$ be a coadjoint orbit of $G$. Let $\pro_{\f h}^{\f g}: \f g^* \rightarrow \f h^*$ be defined as 
$$(\pro_{\f h}^{\f g} \phi)(h)=\phi(h), \qquad (\forall \,\,\, \phi \in \f g^*, \,\,\, h \in \f h).$$
Clearly the set $\pro_{\f h}^{\f g} \mathcal O(\lambda)$ is $H$-invariant under the coadjoint action. Hence 
$\pro_{\f h}^{\f g} \mc O(\lambda)$ is a union of  coadjoint orbits of $H$. For each $\lambda \in \mathfrak g^*//G$, we define $\pro_{\f h}^{\f g} (\lambda)$ to be the subset of $\f h^*//H$, 
$$\{ \mu \in \f h^*//H: \mc O( \mu) \in \pro_{\f h}^{\f g} (\mc O(\lambda))\}. $$
Notice $\pro_{\f h}^{\f g} \mc O(\lambda)$ is a subset of $\f h^*$ and  $\pro_{\f h}^{\f g} \lambda$ is a subset of $\f h^*//H$. \\
\\
It may have been  more reasonable to use ${\rm Res}_{\f h}^{\f g}$, instead of ${\pro}_{\f h}^{\f g}$.  However, in most applications in the literature, $\f g$ and $\f h$ are both reductive. In this situation,  by identifying $\f g^*$ with $\f g$ and $\f h^*$ with $\f h$, $\pro_{\f h}^{\f g}$ is indeed a projection from $\f g$ to $\f h$. 
Since we focus on the case both $G$ and $H$ are reductive, we retain the notation $\pro_{\f h}^{\f g}$. One advantage we gain is that $\f g$ and $\f h$ can both be treated as linear Lie algebras. \\
\\
The most well-known case of orbital projection is the following
\begin{theorem}[Cauchy Interlacing Relation]\label{cauchy}
Let $i \f u(n)$ be the space of $n \times n$ Hermitian matrices. Define $\pro_{i \f u(n-1)}^{i \f u(n)}$ to be the projection of a Hermitian matrix to its upper left $(n-1) \times (n-1)$ submatrix. Then $\mu \in \pro_{i {\f u}(n-1)}^{i {\f u}(n)}(\lambda)$ if and only if $(\lambda, \mu)$ satisfies the Cauchy interlacing relations:
$$\lambda_1 \geq \mu_1 \geq \lambda_2 \geq \mu_2 \geq \ldots \geq \lambda_{n-1} \geq \mu_{n-1} \geq \lambda_n.$$
\end{theorem}
The coadjoint orbits become homogeneous symplectic manifolds when they are equipped with the canonical Kirillov-Kostant-Souriau symplectic form. The orbital projection $\pro_{\f h}^{\f g} \mc O(\lambda)$ can then be interpreted as the moment map of the Hamiltonian $H$-action on $\mc O(\lambda)$. The case that both $G$ and $H$ are compact was intensively studied in the past. In this case, the $\pro_{\f h}^{\f g}(\lambda)$ is a convex polytope, often known as the Kirwan polytope,  denoted by $\Delta_{H}(\lambda)$. We shall refer the reader to \cite{gs} for references and historical account. \\
\\
In our paper, we will focus on the case where both  $G$ and $H$
are noncompact. 
When $G$ is semisimple and noncompact, $\pro_{\f h}^{\f g}(\mc O(\lambda))$ is a lot more complex. First of all, there are often  nonconjugate Cartan subalgebras and each Cartan subalgebra  yields a class of semisimple adjoint orbits. Each class of semisimple adjoint orbit must be treated differently. Secondly, semisimple adjoint orbits do not exhaust all adjoint orbit and there are nilpotent orbits which behave quite differently. Thirdly to treat all adjoint orbits, an induction process will be needed and it will involve the nilpotent orbits of certain smaller subalgebra. 
Generally speaking, $\pro_{\f h}^{\f g}(\lambda) $ will not be convex unless $H$ is compact.
\section{Branching Laws}
Let $G$ be a Lie group and $\pi$ be a unitary representation of $G$.
Let $H$ be a Lie subgroup of type I and $\hat H$ the unitary dual of $H$ (\cite{dix}). Then $\pi|_{H}$ decomposes into a direct integral of irreducible unitary representations of $H$ with multiplicities
$$\int_{\mu \in \hat{H}} \oplus^{m_{\pi}(\mu)} \mu d_{\pi} (\mu).$$
Here $m_{\pi}(\mu)$ is the multiplicity of $\mu$ in $\pi$, which can assume the value $\infty$. Perhaps, the most important part of this direct integral decomposition is the discrete spectrum, namely, those $\mu \in \hat{H}$ that occur as  subrepresentations of $\pi|_H$. We denote the discrete spectrum by $\pi|_{H}^{dis}$. \\
\\
The fundamental case to study is the decomposition of $\pi|_H$  when $\pi$ is irreducible. For irreducible $\pi$, the direct integral decomposition of $\pi|_H$ is often called a branching law. One of the most well-known branching laws is the following
\begin{theorem}[Weyl]\label{weyl}
Let $\lambda \in \mathbb Z^n$ be arranged in descending order. Let $\pi_{\lambda}$ be the irreducible unitary representation of $U(n)$ with highest weight $\lambda$. Let $U(n-1)$ be any subgroup of $U(n)$ preserving a nonzero vector in $\mathbb C^n$. Then
$$\pi_{\lambda}|_{U(n-1)} \cong \bigoplus_{{\mu \in \mathbb Z^{n-1},\lambda_1 \geq \mu_1 \geq \lambda_2 \geq \mu_2 \geq \ldots \geq \lambda_{n-1} \geq \mu_{n-1} \geq \lambda_n}}  \tau_{\mu}.$$
Here $\tau_{\mu}$ is an irreducible unitary representation of $U(n-1)$ corresponding to the highest weight $\mu$.
\end{theorem}
We observe here that Cauchy interlacing relation is manifested here in the branching law.\\
\\
Here we list a few important cases of branching laws.
\begin{enumerate}
\item Let $(G,H)$ be a symmetric pair, i.e., the subgroup $H$ is the set of fixed points of a certain involution of $G$. Let $\pi$ be an irreducible unitary representation of $G$. If $\pi$ is an induced representation, then the problem of $\pi|_H$ can be reduced to a direct integral decomposition of a Hilbert space on a homogeneous vector bundle. However, when $\pi$ is a discrete series representation, $\pi|_H$ can be difficult to study. Some important results were obtained by Kobayashi in a series of papers (\cite{Ko1}\cite{Ko2}\cite{Ko3}) regarding the discrete decomposibility of a larger class of unitary representation $A(\f q, \lambda)$.
\item Let $(G,H)$ be both compact. Then every $\pi \in \hat G$ is finite dimensional and is in the discrete series. This case has  been intensely studied. A particular example is $G=U(n) \times U(n)$ and $H=U(n)$ diagonally embedded in $G$ (\cite{kt}). The branching law in this case is exactly the decomposition of tensor product of two irreducible unitary representations, given by Littlewood-Richardson rule.
\item Let $(H_1, H_2)$ be a dual reductive pair in a symplectic group $Sp_{2n}(\mathbb R)$ and $\omega$ the Weil representation of $\widetilde{Sp}_{2n}(\mathbb R)$. Then branching law $\omega|_{\widetilde{H_1 \times H_2}}$ is given by the $L^2$ Howe's correspondence (\cite{howe}). We understand well the cases that $(H_1, H_2)$ are of the same size (\cite{p} \cite{ab}) or $(H_1, H_2)$  in the stable range (\cite{li2}) or one of $(H_1, H_2)$ is compact (\cite{kave}). There are many other cases that we do not have a complete description of the branching law (\cite{lid}).
\item There are also the problems of studying the multiplicities in the branching laws. A great deal of effort was made to understand the multiplicity one branching laws, both in compact cases and noncompact cases. For the compact case, see \cite{gw} and the references therein.
\end{enumerate} 
The literature for the branching laws  is vast.  It is impossible to include all the references here. The main point is that many of these branching problems have  "counterpart" problems for coadjoint orbits, as we shall discuss in the next section.
\section{Orbit Method}
From Theorem \ref{cauchy} and Theorem \ref{weyl}, we already see that orbit projection and the  branching law are related. In fact, this kind of relation is exactly suggested by  the orbit method. The orbit method, pioneered by  Kirillov and Kostant, is a method that produces a correspondence between coadjoint orbits and irreducible unitary representations, by which both the structure of coadjoint orbits and irreducible unitary representations can be better understood. \\
\\
In the early 1960's, Kirillov proved that there is a one-to-one correspondence between irreducible unitary representations of a simply connected nilpotent group $N$ and the coadjoint orbits in $\f n^*$ (\cite{kir0}). Later Auslander and Kostant extended this correspondence to type I solvable groups (\cite{ak}).  Auslander and Kostant proved that the correspondence between coadjoint orbits and irreducible unitary representations still holds if one imposes an admissible condition on the coadjoint orbits. For $G$ noncompact and semisimple, the correspondence between $\hat G$ and $\f g^*//G$ runs into several problems. For example, not every coadjoint orbit will have a corresponding irreducible unitary representation. Hence there is no one-to-one correspondence between $\f g^*//G$ and $\hat G$. In \cite{du}, Duflo defined an admissible condition on coadjoint orbits for all semisimple Lie groups. This produces an remedy---one should use the admissible orbit datum, instead of coadjoint orbit alone, to associate to irreducible unitary representations (\cite{du}\cite{vogan}). Since then, a great deal of progress has been made to establish a correspondence between the unitary dual and admissible coadjoint orbit datum (\cite{vogan}).\\
\\
 The most well-understood part of the orbit method is the correspondence between the tempered dual and  admissible coadjoint orbits (\cite{du82}).  
In the case that $G$ is real reductive, a theorem due to Harish-Chandra and Rossmann  provides a good foundation to discuss irreducible tempered  representations and  coadjoint orbits (\cite{hc, HCi} \cite{R1} \cite{vec}). We remark that the correspondence between tempered representations and admissible coadjoint orbits have been extended to all type I real Lie groups by Duflo and his school, namely the decomposition of $L^2(G)$ can be parametrized by admissible coadjoint orbits in the sense of Duflo (\cite{du82} \cite{ver}).
\begin{theorem}[Rossmann]
Let $G$ be a real reductive group of Harish-Chandra class. Let $\pi$ be an irreducible tempered representation of $G$ with regular infinitesimal character (\cite{kn}). Then there exists a unique coadjoint orbit $\mc O(\lambda)$ such that $\Theta_{\pi}(g)$, the Harish-Chandra character of $\pi$   satisfies
$$\Theta_{\pi}(\exp x)=C_{\pi} p^{-1}(x) \mathcal F(\mu(\mathcal O(\lambda)))(x) , \qquad (x \in \Omega)$$
as distributions on $\f g$.
Here $C_{\pi}$ is a constant; $\mu(\mc O(\lambda))$ is the canonical invariant measure on $\mc O(\lambda)$; $p^{-1}(x)$ is a factor related to $\exp: \f g \rightarrow G$; $\Omega$ is an open neighborhood of zero in $\f g$  and 
$\mc F$ is the Fourier transform from tempered distributions $\mc S^{\prime}(\f g^*)$ to $\mc S^{\prime}(\f g)$.
\end{theorem}
Essentially, the pull back of the Harish-Chandra character $\exp^*(\Theta_{\pi})$, as a distribution on $\Omega$, can be identified locally with the Fourier transform of the invariant measure $d \mc O_{\lambda}$ multiplied by $C_{\pi} p^{-1}(x)$.
Alternatively, for any $\phi$ a smooth function with compact support in a small open neighborhood of the identity in $G$,
$$ \int \phi(\exp x) \Theta_{\pi}(\exp x) d x=C_{\pi} \int_{\mc O_{\lambda}} \mc F( p^{-1}(x) \phi(\exp x)) (\xi) d \xi.$$

Rossmann's character formula only relates the Harish-Chandra character to the Fourier transform of an admissible orbit in a small neighborhood of the identity.  M. Vergne established the following elegant result, known as the Poisson-Plancherel formula for all semisimple Lie groups (\cite{ve}), in which the complete Fourier transforms of  generic admissible orbits all appear.
\begin{theorem}[Vergne] Let $G$ be a semisimple Lie group. Let
$$E(G)=\{ x \in \f g \mid \exp x=e \}.$$
Then there exist a $G$-invariant measure $d_E$ on $E(G)$ and a $G$-invariant measure $d_v$ on $\f g^*$ such that for any smooth function $\phi$ with compact support
$$\int_{E(G)} \phi(x) d_E(x)= \int_{\f g^*} \mc F(p(x) \phi( x))(\xi) d_v(\xi).$$ 
\end{theorem}
Here $E(G)$ is a (lattice-like) subset of elliptic elements in $\f g$ and $d_v(g)$ can be regarded as a Lie algebra analog of the Plancherel measure
on $\hat G$. Recall that the classical Plancherel formula related
$\phi(e)$ to an integral of $Tr(\pi(\phi))$ with respect to the Plancherel measure on $\pi \in \hat{G}$. \\
\\
%%In any case, what the orbit method suggests is that the unitary dual $\hat G$ is closely related to the coadjoint orbit $\f g^*//G$ and understanding coadjoint orbits can help us understand unitary representations. As we have seen so far, to a large extent, the orbit method is well-established for regular semisimple coadjoint orbits and tempered representations.\\
There are many other important results concerning coadjoint orbits and irreducible unitary representations, which are beyond the scope of this paper.
What is important for  the purpose of this paper, is that, associated to  each discrete series representation of a semisimple group, there is a unique elliptic coadjoint orbit, generated by its Harish-Chandra parameter, from which many analytic properties of the representation can be obtained. \\
\\
To understand what the orbit method says about branching law, let us recall the theorem of Kirillov regarding branching laws of nilpotent groups (page 81, \cite{kir0}).
\begin{theorem}[Kirillov] Let $N$ be a simply connected nilpotent Lie group and $H$ a closed connected Lie subgroup. For each $\lambda \in n^*//N$, let $\pi_{\lambda}$  be the irreducible unitary representation constructed by polarization of the coadjoint orbit $\mc O(\lambda)$. Then 
$$\pi_{\lambda}|_H =\int_{\mu \in {{\pro}_{\f h}^{\f n} (\lambda)}} \oplus^{m(\mu)}\tau_{\mu} d_{\lambda} (\mu),$$
with $\tau_{\mu}$ the irreducible unitary representation of $H$ corresponding to $\mu \in \f h^*//H$. 
\end{theorem}
This theorem can be characterized by the following commutative diagram
\begin{equation}
\begin{CD}
\mc O(\lambda) \subset \f n^*  @>{}>> \pi_{\lambda} \\
@VV{\pro_{\f h}^{\f n}}V    @VV{Res|_{H}}V \\
\sqcup_{\mu \in \pro_{\f h}^{\f n}(\lambda)} \mc O(\mu) @>{}>> \int_{\mu \in \pro_{\f h}^{\f n}(\lambda)} \oplus^{m(\mu)} \tau_{\mu} d_{\lambda} (\mu)
\end{CD}
\end{equation}
The remaining question is whether Kirillov's theorem can be extended to other groups. \\
\\
When $G$ is compact, $H$ is a closed Lie subgroup, $\pi$ an irreducible unitary representation, $\pi|_H$ decomposes discretely into a direct sum of irreducible representations of $H$ with multiplicities:
$$\pi|_H=\oplus_{\sigma \in \hat H} \sigma \otimes M(\pi, \sigma).$$
Here $M(\pi, \sigma)$ records the multiplicity. On the geometric side, there is the moment map, which can be identified with $\pro_{\f h}^{\f g}$, maps any coadjoint orbit $\mc O$ into $\f h^*$. The famous Guillemin-Sternberg conjecture, proved by Meinrenken, established a precise relation between $M(\pi, \sigma)$ and the geometry of the preimage of $\mc O(\sigma)$ (\cite{gst} \cite{heck} \cite{me, me1}). Here the correspondence between coadjoint orbits and unitary dual is readily provided by the Borel-Weil theorem. \\
%%A lot of work has been done for $G$ is compact, starting with Heckman's paper  (\cite{heck}) and the Guillemin-Sternberg conjecture ( \cite{gst}). \cite{kir})). \\
\\
When $G$ is semisimple noncompact and $H$ is compact, $\pi$ a tempered representation of $G$, Guillemin-Sternberg conjecture remains to be true with some modifications (\cite{pa} \cite{dve}).
As pointed out by us by Prof. Vergne, what holds true for $H$ compact, still applies to the case $H$ noncompact under the  assumption that the projection $\pro_{\f k_{H}}^{\f g}: \mathcal O \rightarrow \f k_{H}^*$ is proper where $K_H$ is a maximal compact subgroup of $H$. In such case, the representation corresponding to $\mc O$ should decompose discretely under $K_H$, therefore decompose discretely under $H$ (\cite{pa1}).\\
\\
 Recently, Liu, Oshima and Yu established a branching law for
$G=Spin(n,1)$ and $H$ the minimal parabolic subgroup. They proved that for any irreducible unitary representation $\pi$ of $G$, 
$\pi|_H$ decomposes into a finite direct sum.  In this case, they  confirmed a conjecture of Duflo , relating the branching law of the discrete series to $\pro_{\f h}^{\f g}$ of
the elliptic orbits (\cite{loy}).\\
\\
For $(G, H)$  noncompact, the relation between the branching laws and orbit projection is much more delicate.   Nevertheless, many techniques for compact $H$ still make sense. For instance, when $\pro_{\f h}^{\f g}: \mc O \rightarrow \f h^*$ is proper, the push-forward of the Louiville measure, often known as the Duistermaat–Heckman measure, still exists (\cite{dh}, \cite{dhv}). In some sense, the Duistermaat–Heckman measure is still closely related to the multiplicity space $M(\pi|_H, \sigma)$ even though a few subtleties, like $\rho$-shift, continuous spectrum and nontemperedness, may complicate the relation. In any case, as we shall see later, the study of orbit projection problem  can often provide new insights and perspectives to the study of branching problem and vice versa. 

%%The case that $\pi$ is a discrete series representation seems to be the most interesting case with least complications 
\section{Discrete Series and Elliptic Orbits}
Let $G$ be a semisimple Lie group. Let $G \times G$ act on $L^2(G)$ from the left and right. Then $L^2(G)$ becomes a unitary representation of $G \times G$. The decomposition of $L^2(G)$ into the direct integral of irreducible unitary representations of $G \times G$ is known as the Plancherel formula.
In the 1960's, Harish-Chandra successfully carried out the determination of the Plancherel formula. One critical step is the determination of the discrete spectrum $L^2(G)^{dis}$. Harish-Chandra proved that
discrete spectrum exists if and only if the group $G$ has a compact Cartan subgroup. He also parametrized the discrete spectrum of $L^2(G)$.  We now state Harish-Chandra's parametrization of discrete series (\cite{HC5} \cite{HCd}).
\begin{theorem}[Harish-Chandra] Let $G$ be a connected reductive Lie group of Harish-Chandra class. Let $K$ be a maximal compact subgroup of $G$. Then $G$ has discrete series representation if and only if $rank(G)=rank(K)$. Let $T$ be a maximal torus of $K$. Let $L$ be the integral lattice dual to $T$ and $\rho$ be half of the sum of positive roots of $\f g$ with respect to $\f t$.  Let $W(K, \f t)$ be the Weyl group of $K$ with respect to $\f t$. Let $(L+\rho)^{\prime}$ be the regular elements in $L+\rho$ with respect to the action of $W(\f g_{\mathbb C}, \f t)$.  Then discrete series of $G$ is in one-to-one correspondence with $(L+ \rho)^{\prime} // W(K, \f t)$, the $W(K, \f t)$ orbits of $(L+\rho)^{\prime}$.
\end{theorem} 
The element  $\lambda \in (L+ \rho)^{\prime} // W(K, \f t)$, or simply $\lambda \in (L+ \rho)^{\prime}$ is called the Harish-Chandra parameter of the corresponding discrete series $D_{\lambda}$. \\
\\
Recall that the semisimple adjoint orbits can be represented by elements in the Cartan subalgebras of $\f g$. Hence they are parametrized by the related Weyl group orbit in the related Cartan subalgebra. The elliptic adjoint orbits are simply orbits generated by elements in the compact Cartan subalgebra. Identify $\f g$ with $\f g^*$ under the Killing form.  $L+\rho$, as a subset of $\f t^*$,  becomes a subset in $\f g^*$. The orbits in
$\f t^*//W(K)$ are then in one-to-one correspondence with the elliptic coadjoint orbits in $\f g^*$.  Now we can attach the discrete series $D_{\lambda}$ to the elliptic coadjoint orbit $\mc O(\lambda)$.  \\
\\
Despite Harish-Chandra's great success of classifying discrete series representation, the structure of the discrete series representation is still not well-understood. In particular, we do not have a good understanding of $K$-types of discrete series representations and their multiplicities. The branching laws of discrete series $\pi|_K$ are provided by the Blattner's formula (\cite{HS}). However, the coefficients in Blattner's formula are difficult to compute and it is not easy to know which ones are nonzero. The related problem on the orbit side,  $\pro_{\f k}^{\f g}$ of the elliptic orbit $\mc O(\lambda)$, is equally difficult. For some recent progresses see \cite{pa} \cite{ha} \cite{BV} and Atlas of Lie groups project.\\
\\
As to restrictions to noncompact subgroups, even though discrete series representation can be constructed geometrically as cohomology classes (\cite{schmid}), the restrictions of cohomolgy classes to submanifolds may not even make sense. There is no general construction of discrete series by which branching laws can be easily derived. One exception here is the holomorphic discrete series. A range of branching laws can be established for  discrete series by Duflo-Vargas (\cite{DV}), Vargas (\cite{var,var1}), Paradan (\cite{pa}), Liu-Oshima-Yu(\cite{loy}), Baldoni-Vergne (\cite{BV1}) among others.  Most efforts have been focused on the case that $\pi$ is $H$-admissible in the sense of Kobayashi.  To this end,  Duflo made a precise conjecture relating the $H$-admissibility to the properness of the moment map
$m: \mc O(\lambda) \rightarrow \f h^*$. See \cite{loy} for the precise statement of Duflo's conjecture.\\
\\
%%Besides using cohomolgy to construct discrete series,  theta correspondence  can also be used to construct discrete series.  There are very few other tools available to study the branching laws of discrete series representations.
%%Get back to the Harish-Chandra parameters. We shall now discuss semisimple adjoint orbits, in particular the elliptic orbits, with the understanding that  the  coadjoint orbits are identified with adjoint orbits once we fix an invariant bilinear form. The reason we want to use adjoint orbits is that the computation on projection of adjoint orbits can be carried out entirely using linear algebra.  \\

\section{The Gan-Gross-Prasad Branching problem for unitary groups}
Let us now consider the discrete series of $G=U(p,q)$.  The compact Cartan subalgebra can be identified with $i \mathbb R^p \times i \mathbb R^q$. In linear algebra terms, an elliptic element in $\f u(p,q)$ is a diagonalizable matrix in $\f u(p,q)$ with purely imaginary eigenvalues
$$(i\chi_1, i\chi_2, \ldots, i\chi_{p+q}).$$
The regular elements are just those $i \chi$ with distinct $\chi_i$. Suppose that $\chi_i$'s are all distinct.
The defining Hermitian form of signature $(p,q)$ restricted onto each eigenspace $E( i \chi_k)$ is either positive definite or negative definite. If we attach a $+$ sign ($+1$) or $-$ sign ($-1$) to $i \lambda_k$,  we obtain a sequence of sign in $\{ \pm 1\}^{p+q}$ with $p$ $+1$'s and $q$ $-1$'s. We denote it by $z$. Now we see that the regular elements in a compact Cartan subalgebra can be parametrized by $(\chi, z)$ with $\chi \in \mathbb R^{p+q}$ and $z \in \{\pm 1\}^{p+q}$ and $z$ has signature $(p,q)$.  The Weyl group $W(K, \f t)$ acts on $(i \chi, z)$ by permuting the $(i \chi, z)$ simultaneously and preserving $z$. In other words, $W(K, \f t)$ permutes those $i \chi_i$ with the same sign $z_i$. For convenience, one may think of $(i \chi, z)$ as signed eigenvalues, with each eigenvalue $i \chi_i$ a sign $z_i$ attached to it. Now an regular elliptic orbit $\mc O(i \chi, z)$ is simply the conjugacy class of matrices in $\f u(p,q)$ with the  signed eigenvalue $(i \chi, z)$.\\
\\
We make a quick remark here. If $\chi$ is not regular,  $E(i \lambda_k)$ may be more than one dimensional. In this case, we will have $\lambda_k$ appear $\dim(E(i \lambda_k))$ times. Suppose that the defining Hermitian form restricted onto $E(i \lambda_k)$, which is necessarily nondegenerate, has signature $(r,s)$. Then we assign $r$ $+1$'s and $s$ $-1$'s to $\lambda_k$'s. There are $\frac{(r+s)!}{r ! s !}$ ways to do this. Hence $(i \chi, z)$ will no longer be unique. In any case, elements in a compact Cartan subalgebra can still be represented by its signed eigenvalues $(i \chi, z)$. There is no ambiguity what $\mc O(i \chi, z)$ is even when $(i \chi,z)$ is not regular. 

\begin{definition}[\cite{gp} \cite{ggp}] We say that two signed elliptic element  $(i\chi, z)$ and $(i\eta, t)$, of $U(p,q)$ and $U(p-1,q)$ respectively, satisfy the Gan-Gross-Prasad interlacing relation,  if one can line up $\chi$ and $\eta$ in the descending ordering such that the corresponding sequence of signs from $z$ and $t$ only has the following eight adjacent pairs
$$(\oplus +), (+ \oplus), (- \ominus), (\ominus -), (+-), (-+), (\oplus \ominus), (\ominus \oplus).$$
Here $\oplus$ and $\ominus$ represent $+1$ and $-1$ in $t$, and $+$ and $-$ represent $+1$ and $-1$ in $z$.
We call such a sign sequence the (interlacing) sign pattern of $(\chi, z)$ and $(\eta, t)$. Clearly, when there is neither $-$ nor $\ominus$, this interlacing relation is exactly the classical  Cauchy interlacing relation.
\end{definition}
Now we have the following theorem regarding the projection of elliptic orbits in $\f u(p,q)$.
\begin{theorem}
Let $(i \chi, z)$ be a regular elliptic element in a compact Cartan subalgebra in $\f u(p,q)$. Suppose that $q \geq 1$. Then
an elliptic orbit $\mc O(i \eta, t)$ appears in $\pro_{\f u(p,q-1)}^{\f u(p,q)} \mc O(i \chi, z)$ if and only if $(i\chi, z)$ and $(i\eta, t)$ satisfy the Gan-Gross-Prasad interlacing relation.
\end{theorem}
The proof of this theorem is based on linear algebra. We omit it here. \\
\\
Let us look at the elliptic orbit corresponding to holomorphic discrete series. Consider $(i \chi, z)$ with
$$\chi_1 > \chi_2 \ldots > \chi_{p+q},$$
$$z_1=z_2= \ldots z_p=1, \qquad z_{p+1}=z_{p+2}= \ldots=z_{p+q}=-1.$$
We have the sign pattern of $(i\chi, z)$:
$$\overbrace{+ \, + \, \ldots \, +}^p \overbrace{ \, - \, -\, \ldots -}^q.$$
We now insert $(i\eta, t)$ into $(i \chi, z)$ in descending order. The only sign pattern allowed by GGP interlacing relation is
$$\oplus + \oplus + \oplus + \ldots  \oplus + - \ominus - \ominus - \ldots  \ominus -.$$
Hence the elliptic element appears in the projection must have the sign pattern
$$\overbrace{\oplus \oplus \ldots \oplus}^{p} \overbrace{\ominus \ominus \ldots \ominus}^{q-1}.$$
We see that these elliptic elements again correspond to holomorphic discrete series. In addition, we can show that these elliptic elements exhaust all regular semisimple orbits in $\pro_{\f u(p,q-1)}^{\f u(p,q)}(\mc O(\xi, z))$.
 
%%There are generally other semisimple orbits in 
%%$\pro_{\f u(p,q-1)}^{\f u(p,q)} \mc O(\i \chi, z)$.  Very briefly, for %%each $\oplus \ominus$ or $\ominus \oplus$ pair, one can substitute
%%the corresponding pair $(i \chi_{k}, i \chi_{k+1})$ by a semisimple %%%%%part with eigenvalue $( i \chi_{k}^*+ \mu, i \chi_{k}^*- \mu)$. The %%%resulting semisimple orbit remains in the projection of $\mc O(\i \chi, z)$.

A similar statement holds on the branching law side, as conjectured by Gan, Gross and Prasad (\cite{gp} \cite{ggp}). Let $(i\chi, z)$ be a Harish-Chandra parameter for a discrete series representation $D(i \chi, z)$ for $U(p,q)$. Here $\chi \in \mathbb R^{p+q}$ is a sequence of distinct integers if $p+q$ is odd, of half integers if $p+q$ is even. Let $D(i \eta, t)$ be a discrete series representation of $U(p-1, q)$.
 \begin{theorem}[\cite{ggphe}] Suppose $q \geq 1$. The discrete spectrum 
$$D(i\chi, z)|_{U(p,q-1)}^{dis} = \hat{\oplus}_{(\eta, t)} D(i\eta, t)$$
where the direct sum is taking over all those Harish-Chandra parameters $(i\eta, t)$ such that $(i\chi, z)$ and $(i\eta, t)$ satisfy the GGP interlacing relation.
 \end{theorem}
 The main idea of the proof is to relate the branching law of $D(i\chi, z)$ to the branching laws of the Weil representation $\omega$ restricted to a noncompact
  dual pair $(U(n), U(n+1))$.  
 One critical ingredient is the determination of discrete spectrum of $\omega$, due to Jian-Shu Li (\cite{lid}). This allows us to apply Howe's correspondences in a see-saw fashion to obtain the discrete spectrum of $D(i\chi, z)|_{U(p,q-1)}$ inductively. \\
 \\ 
As we can see, in the case of $(U(p,q), U(p, q-1))$, the projection of elliptic coadjoint orbits $\pro_{\f h}^{\f g}$ does tell correctly which discrete series representations occurs in $D_{H}$. However, one cannot use this fact to prove the branching law. The orbit method is a good way to suggest the branching laws. But one can hardly use it to prove the branching laws when $G$ is semisimple (\cite{kir}). Kirillov's branching law is an exception since all irreducible unitary representations can be obtained through polarization. 

\section{Wave Front set}
The projection of coadjoint orbits $\pro_{\f h}^{\f g}$, in the most general cases, will not tell exactly which representations occurs in $\pi|_{H}$. In addition, the correspondence between coadjoint orbits and unitary dual may not be available. One remedy is to study the wave front set.
Let $\pi$ be a unitary representation of $G$, not necessarily irreducible. The wave front set $WF(\pi)$, defined by  Howe, is the union of the wave front sets of all matrix coefficients of $\pi$ (\cite{How}). It is a conic subset of the cotangent bundle $T^* G$. Due to the action of $G$, it is enough to consider $T^*_e G$ which can then be identified with $\f g^*$. Without loss of generality, we assume $WF(\pi)$ is in $\f g^*$. Due to the adjoint action of $G$, $WF(\pi)$ is invariant under the coadjoint action of $G$.  \\
\\
When $\pi$ is an irreducible unitary representation of a semisimple Lie group, $WF(\pi)$  lies in the zero set of all homogeneous $G$-invariant functions on $\f g^*$. Hence $WF(\pi)$ lies in $\mc N(\f g^*)$, the nilcone of $\f g^*$. In particular, $WF(\pi)$ can be defined as  the wave front set of the Harish-Chandra character $\Theta_{\pi}(g)$. It is closed related to other invariants of $\pi$, the asymptotic cycle, characteristic cycle, and associated variety (\cite{bv0} \cite{SV} \cite{V}).  Since a series of representations may share the same wave front set, one usually does not gain precise information about the branching laws. The issue is less about $G$,  since many irreducible representations of $G$ with the same wave front set also have the similar branching laws. It is more about the group $H$. More precisely, let us consider the following diagram
\begin{equation}
\begin{CD}
\pi \subset \hat{G}  @>{WF}>> WF(\pi) \\
@VV{Res}V    @VV{\pro}V \\
\pi|_H=\int_{\mu \in \hat H} \oplus^{m(\mu)} {\mu} d_{\lambda} (\mu) @>{???}>> \pro_{\f h}^{\f n} WF(\pi)
\end{CD}
\end{equation}
One important result, due to Howe, asserts that
$$WF(\pi|_H) \supseteq \pro_{\f h}^{\f g} WF(\pi).$$
But we do not know whether $WF(\pi|_H) = \pro_{\f h}^{\f g} WF(\pi)$. Even if we assume the equality, we still cannot read $\pi|_H$ off the
$WF(\pi|_H)$. 
Nevertheless, once we know $\pro_{\f h}^{\f g} WF(\pi)$, we gain some information about  $\pi|_H$. There are many important results regarding branching laws based on wave front sets.
\begin{theorem}[Kobayashi, Cor. 3.4 \cite{Ko3}]
Let $\pi$ be an irreducible unitary representation of a reductive group $G$. Let $H$ be a reductive subgroup. Suppose that $\pi|_H$ is infinitesimally discretely decomposable, then
$$\pro_{\f h}^{\f g}(WF(\pi)) \subseteq \mc N(\f h^*).$$
Here reductive subgroup means $H$ is reductive in $G$. Infinitesimally discretely decomposable means that the underlying $(\f g, K)$ module can be decomposed as a direct sum of irreducible $(\f h, K_H)$ submodules.
\end{theorem}
Clearly, this theorem gives us a criterion for $\pi|_H$ not infinitesimally discretely decomposable. \\
\\
For the branching law $\omega|_{\tilde G}$ in the dual reductive setting, we have
\begin{theorem}[Przebinda]
Let $(G, G^{\prime})$ be a dual pair $Sp_{2N}(\mathbb R)$. Let $\omega$ be the Weil representation of $\widetilde{Sp}_{2N}(\mathbb R)$. If $\pi$ is in the discrete spectrum of $\omega|_{\tilde G}$ then 
$$WF(\pi) \subseteq \pro_{\f g}^{\f{sp}_{2N}(\mathbb R)} WF(\omega).$$
Here the $WF(\omega)$ is the minimal nilpotent orbit of $\f{sp}_{2N}(\mathbb R)$ consisting of all rank 1 and 0 matrices in $\f{sp}_{2N}(\mathbb R)$.
\end{theorem}
In fact, Przebinda proved this theorem for all $\pi \in \mc R(\tilde G, \omega)$( Cor. 2.8, page 557 \cite{pr}). This theorem can be applied to detect those representations that are not in $\mc R(\tilde G, \omega)$ (\cite{hhequan}). Important results concerning $\pro_{\f g}^{\f{sp}_{2N}(\mathbb R)} WF(\omega)$ can be found in \cite{dkp}.\\
\\
Finally, we shall mention a result due to Harris, Olafsson and myself that describes the wave front set for any unitary representation weakly contained in the $L^2$ space of a reductive group of Harish-Chandra class (\cite{hho}).
Let $G$ be a reductive Lie group of Harish-Chandra class. Let $\widehat{G}_{\text{temp}}$ be the tempered dual (\cite{KZ}), the part of unitary dual appearing in the Plancherel formula of $L^2(G)$. To each irreducible tempered representation $\sigma$ of $G$, as we have seen for the discrete series, Duflo and Rossmann associated a finite union of coadjoint orbits $\mathcal{O}_{\sigma}\subset \mathfrak{g}^*$ ( \cite{Du},\cite{R1},\cite{R2}). In the generic case, $\mathcal{O}_{\sigma}$ is a single coadjoint orbit.\\
\\
For each $\pi$ weakly contained in the regular representation, we define the orbital support of $\pi$ by 
$$\mathcal{O}\operatorname{-}\operatorname{supp}\pi=\bigcup_{\sigma\in \operatorname{supp} \pi}\mathcal{O}_{\sigma}.$$
Here $\operatorname{supp} \pi$ is a closed subset of $\Hat G_{\text{temp}}$.

\begin{theorem}  If $G$ is a noncompact reductive Lie group of Harish-Chandra class and $\pi$ is weakly contained in the regular representation of $G$, then
$$WF(\pi)=AC(\mathcal{O}\operatorname{-}\operatorname{supp} \pi).$$
where $AC(S)$ for any $S$ in a linear space $V$
is defined to be
$$AC(S)=\{\xi \in V \mid  \Gamma\ \text{an\ open\ cone\ containing}\
\xi \Rightarrow \Gamma\cap S\ \text{is\ unbounded} \} \cup \{0\}.$$
\end{theorem}

When $G$ is compact and connected, this result is known. See Cor 5.10 of \cite{KV}, Proposition 2.3 of \cite{How} and \cite{HHK}. 
\begin{corollary} Let $\pi$ be a unitary representation of a reductive Lie group $G$ in the Harish-Chandra class. Suppose that $WF(\pi)$ contains a regular elliptic element in $\f g^*$. Then there are infinitely many discrete series representations in the discrete spectrum of $\pi$.
\end{corollary}

Combined with Howe's result that $WF(\pi|_H) \supseteq \pro_{\f h}^{\f g} WF(\pi)$, we have
\begin{corollary} Let $G$ be a  Lie group and $H$ be a Lie subgroup that is reductive of the Harish-Chandra class. Let $\pi$ be an irreducible unitary representation of $G$ such that $\pi|_H$ is weakly contained in $L^2(H)$.  Suppose that $\pro_{\f h}^{\f g} WF(\pi)$  contains a regular elliptic element in $\f h^*$. Then there are infinitely many discrete series representations of $H$ in the discrete spectrum of $\pi|_H$.
\end{corollary}
\section*{acknowlegement} 
I would like to thank Prof. R. Howe for getting me interested in discrete series, Prof. M. Vergne for her very helpful comments and her interests in this work, and J. Yu for bring their work (\cite{loy}) to my attention. I would also like to thank the editors for allowing me revising Section 3 after the original manuscript was accepted.

\bibliographystyle{amsplain}

\end{document}